\documentclass[11pt]{amsart}
\usepackage[colorlinks=true, pdfstartview=FitV, linkcolor=blue, citecolor=blue, urlcolor=blue, breaklinks=true]{hyperref}
\usepackage{amsmath,amsfonts,amssymb,amsthm,amscd,comment,euscript}

\usepackage[usenames]{color}
\usepackage{hyperref}
\usepackage{amsmath, amscd}
\usepackage{tikz}
\usepackage{tikz}
\usepackage[all]{xy}
\usepackage{tikz}
\usepackage{array}
\usepackage{bookmark}
\usepackage{graphicx}

\newcounter{cnt}
 \makeatletter
\def\mydggeometry{\makeatletter\dg@YGRID=1\dg@XGRID=20\unitlength=0.003pt\makeatother}
\makeatother \theoremstyle{remark}
\numberwithin{equation}{section}

\let\bwdg\bigwedge
\def\bigwedge{{\textstyle\bwdg}}
\theoremstyle{definition}

\theoremstyle{definition}
\newtheorem{theom}{Theorem}

\newtheorem*{theoA}{Theorem}

\newtheorem*{coro}{Corollary}

\newtheorem*{remi}{Remark}

\newcommand\al{\alpha}

\newcommand\Lg{\mathfrak{g}}
\newcommand\Lh{\mathfrak{h}}

\newcommand\Ln{\mathfrak{n}}

\newcommand{\om}{{\omega}}
\newcommand{\lam}{{\lambda}}

\newcommand{\lie}{\mathfrak}
\newcommand{\mindeg}{\operatorname{mindeg}}

\newcommand{\chara}{\operatorname{char}}

\begin{document}

\title[The PBW degree of simple modules]{The degree of the Hilbert-Poincar\'e polynomial of PBW-graded modules}
\author[Backhaus, Bossinger, Desczyk, Fourier]{Teodor Backhaus and Lara Bossinger and Christian Desczyk and Ghislain Fourier}
\address{\newline
Mathematisches Institut, Universit\"at zu K\"oln, Germany}
\email{tbackha@math.uni-koeln.de}
\email{lbossing@math.uni-koeln.de}
\email{cdesczyk@math.uni-koeln.de}
\address{\newline
Mathematisches Institut, Universit\"at Bonn, Germany\newline
School of Mathematics and Statistics, University of Glasgow, UK}
\email{ghislain.fourier@glasgow.ac.uk}
\thanks{}

\subjclass[2010]{}
\begin{abstract}
In this note, we study the Hilbert-Poincar\'e polynomials for the associated PBW-graded modules of simple modules for a simple complex Lie algebra. The computation of their degree can be reduced to modules of fundamental highest weight. We provide these degrees explicitly.

\bigskip
Nous \'etudions les polyn\^omes de Hilbert-Poincar\'e pour les modules PBW-gradu\'es associ\'es aux modules simples d'une alg\`ebre de Lie simple complexe. Le calcul de leur degr\'e peut \^etre restreint aux modules de plus haut poids fondamental. Nous donnons une formule explicite pour ces degr\'es.
\end{abstract}
\maketitle \thispagestyle{empty}
%%%%%%%%%%%%%%%%%%%%%%%%%%%%%%%%%%%%%%%%%%%%%%%%%%%%%%%%%%%%%%%%%%%%%%%%%%%%%%%%%%%%%%%%%%%%%%%%%%%%%%%%%%%%%%%%%%%%%%%%%%%%%%%%%%%
%         Introduction
%%%%%%%%%%%%%%%%%%%%%%%%%%%%%%%%%%%%%%%%%%%%%%%%%%%%%%%%%%%%%%%%%%%%%%%%%%%%%%%%%%%%%%%%%%%%%%%%%%%%%%%%%%%%%%%%%%%%%%%%%%%%%%%%%%%
\section{Introduction}
Let $\lie g$ be a simple complex finite-dimensional Lie algebra with triangular decomposition $\lie g = \lie n^+ \oplus \lie h \oplus \lie n^-$. Then the PBW filtration on $U(\lie n^-)$ is given as $U(\Ln^-)_s:= \mathrm{span}\lbrace x_{i_1}\cdot \cdot  \cdot x_{i_l} \mid x_{i_j} \in \Ln^-, l\le s \rbrace$. The associated graded algebra is isomorphic to $S(\lie n^-)$.
Let $V(\lambda)$ be a simple finite-dimensional module of highest weight $\lambda$ and $v_\lambda$ a highest weight vector. Then we have an induced filtration on $V(\lambda) =  U(\Ln^-)v_{\lambda}$, denoted $V(\lambda)_s := U(\Ln^-)_sv_{\lambda}$. The associated graded module $V(\lambda)^a$ is a $S(\lie n^-)$-module generated by $v_\lambda$. \\
These modules have been studied in a series of papers. Monomial bases of the graded modules and the annihilating ideals have been provided for the $\lie{sl}_n, \lie{sp}_n$ \cite{FFoL11a, FFoL11b, FFoL13}, for cominuscule weights and their multiples in other types \cite{BD14}, for certain Demazure modules in the $\lie{sl}_n$-case in \cite{Fou14b, BF14}. In type $G_2$ there is a monomial basis provided by \cite{G11}.\\
The degenerations of the corresponding flag varieties have been studied in \cite{Fei12, FFoL13a, CL14, CLL14}. Further, it turned out (\cite{Fou14}), that these PBW degenerations have an interesting connection to fusion product for current algebras. The study of the characters of PBW-graded modules has been initiated in \cite{CF13, FM14}.\\
In the present paper we will  compute the maximal degree of PBW-graded modules in full generality (for all simple complex Lie algebras), where there have been partial answers in the above series of paper for certain cases. \\
We denote the Hilbert-Poincar\'e series of the PBW-graded module, often referred to as the \textit{$q$-dimension of the module}, by
\[
p_{\lambda}(q) = \sum_{s = 0}^{\infty} \left(\dim V(\lambda)_s / V(\lambda)_{s-1} \right) q^s.
\]
Since $V(\lambda)$ is finite-dimensional, this is obviously a polynomial in $q$. In this note we want to study further properties of this polynomial. We see immediately that the constant term of $p_{\lambda}(q)$ is always $1$ and the linear term is equal to
\[
\dim (\lie n^-) - \dim \operatorname{Ker }\left(\lie n^- \longrightarrow \operatorname{End } (V(\lambda))\right).\]
 Our main goal is to compute the degree of $p_{\lambda}(q)$ and the first step is the following reduction \cite[Theorem 5.3 ii)]{CF13}:
\begin{theoA}
Let $\lambda_1, \ldots, \lambda_s \in P^+$ and set $\lambda = \lambda_1 + \ldots + \lambda_s$. Then
\[
\deg p_{\lambda}(q) = \deg p_{\lambda_1}(q)  + \ldots + \deg p_{\lambda_s}(q).
\]
\end{theoA}
It remains to compute the degree of $p_{\lambda}(q) $ where $\lambda$ is a fundamental weight. We have done this for all fundamental weights of simple complex finite-dimensional Lie algebras:
\begin{theom}\label{B}
The degree of $p_{\omega_i}(q)$ is equal to the label of the $i$-th node in the following diagrams:\\
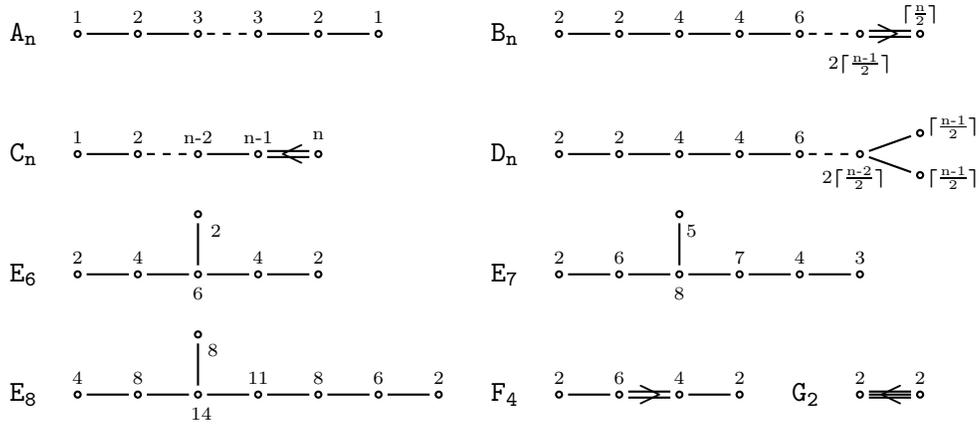
\begin{center}
\begin{tikzpicture}[scale=.4]
    \draw (-21,0) node[anchor=east]  {$\mathtt{A_n}$};
    \foreach \x in {-10,...,-6,-5}
    \draw[thick,xshift=\x cm] (\x cm,0) circle (1 mm);
    \foreach \y in {-10,...,-9,-7,-6}
    \draw[thick,xshift=\y cm] (\y cm,0) ++(.3 cm, 0) -- +(14 mm,0);
    %\draw[thick,dashed] (4 cm,2 cm) circle (1 mm);
    %\draw[thick] (4 cm, 3mm) -- +(0, 1.4 cm);
    %\draw (0,0) node[anchor=north]  {\tiny{1}};
		\draw (-20,1.1) node[anchor=north]  {\tiny{1}};
    %\draw (2,0) node[anchor=north]  {\tiny{2}};
		\draw (-18,1.1) node[anchor=north]  {\tiny{$2$}};
    %\draw (4,0) node[anchor=north]  {\tiny{3}};
		\draw (-16,1.1) node[anchor=north]  {\tiny{$3$}};
    %\draw (6,0) node[anchor=north]  {\tiny{n-2}};
		\draw (-14,1.1) node[anchor=north]  {\tiny{$3$}};
    %\draw (8,0) node[anchor=north]  {\tiny{n-1}};
		\draw (-12,1.1) node[anchor=north]  {\tiny{2}};
		%\draw (10,-0.1) node[anchor=north]  {\tiny{n}};
		\draw (-10,1.1) node[anchor=north]  {\tiny{1}};
    \draw[thick,dashed] (-15.7,0)--(-14.3,0);
    %\draw (4,3) node[anchor=north]  {\tiny{2}};
  %\end{tikzpicture}
%\end{center}
%\bigskip
%\begin{center}
%\begin{tikzpicture}[scale=.4]
    \draw (-5,0) node[anchor=east]  {$\mathtt{B_n}$};
    \foreach \x in {-2,-1,0,...,4}
    \draw[thick,xshift=\x cm] (\x cm,0) circle (1 mm);
    \foreach \y in {-2,-1,0,...,0,1}
    \draw[thick,xshift=\y cm] (\y cm,0) ++(.3 cm, 0) -- +(14 mm,0);
		\draw[thick] (6.3,0.1) -- +(1.4 cm, 0);
		\draw[thick] (6.3,-0.1) -- +(1.4 cm, 0);
		\draw (7.7,0) node[anchor=east] {${{\bf\boldsymbol{>}}}$};
   %\draw (-4,0) node[anchor=north]  {\tiny{1}};
		\draw (-4,1.1) node[anchor=north]  {\tiny{2}};
		%\draw (-2,0) node[anchor=north]  {\tiny{2}};
		\draw (-2,1.1) node[anchor=north]  {\tiny{2}};
    %\draw (0,0) node[anchor=north]  {\tiny{3}};
		\draw (0,1.1) node[anchor=north]  {\tiny{4}};
    %\draw (2,0) node[anchor=north]  {\tiny{4}};
		\draw (2,1.1) node[anchor=north]  {\tiny{4}};
    %\draw (4,0) node[anchor=north]  {\tiny{5}};
		\draw (4,1.1) node[anchor=north]  {\tiny{6}};
    %\draw (6,1.3) node[anchor=north]  {\tiny{n-1}};
		\draw (6,-0.3) node[anchor=north]  {\tiny{$2\boldsymbol\lceil \frac{\textrm{n-1}}{2}\boldsymbol\rceil$}};
    %\draw (8,-0.1) node[anchor=north]  {\tiny{n}};
		 \draw (8,1.3) node[anchor=north]  {\tiny{$\boldsymbol\lceil \frac{\textrm{n}}{2}\boldsymbol\rceil$}};
    \draw[thick,dashed] (4.3,0)--(5.7,0);
    %\draw (4,3) node[anchor=north]  {\tiny{2}};
  %\end{tikzpicture}
%\end{center}
%%\bigskip
%\begin{center}
%\begin{tikzpicture}[scale=.4]
    \draw (-21,-4) node[anchor=east]  {$\mathtt{C_n}$};
    \foreach \x in {-10,...,-6}
    \draw[thick,xshift=\x cm] (\x cm,-4) circle (1 mm);
    \foreach \y in {-10,...,-10,-8}
    \draw[thick,xshift=\y cm] (\y cm,-4) ++(.3 cm, 0) -- +(14 mm,0);
		\draw[thick] (-13.7,-3.9) -- +(1.4 cm, 0);
		\draw[thick] (-13.7,-4.1) -- +(1.4 cm, 0);
		\draw (-12.1,-4) node[anchor=east] {${{\bf\boldsymbol{<}}}$};
    %\draw (0,0) node[anchor=north]  {\tiny{1}};
		\draw (-20,-2.9) node[anchor=north]  {\tiny{1}};
    %\draw (2,0) node[anchor=north]  {\tiny{2}};
		\draw (-18,-2.9) node[anchor=north]  {\tiny{2}};
    %\draw (4,0) node[anchor=north]  {\tiny{n-2}};
		\draw (-16,-2.9) node[anchor=north]  {\tiny{n-2}};
    %\draw (6,0) node[anchor=north]  {\tiny{n-1}};
		\draw (-14,-2.9) node[anchor=north]  {\tiny{n-1}};
    %\draw (8,-0.1) node[anchor=north]  {\tiny{n}};
		\draw (-12,-2.9) node[anchor=north]  {\tiny{n}};
    \draw[thick,dashed] (-17.7,-4)--(-16.3,-4);
    %\draw (4,3) node[anchor=north]  {\tiny{2}};
  %\end{tikzpicture}
%\end{center}
%\bigskip
%\begin{center}
%\begin{tikzpicture}[scale=.4]
    \draw (-5,-4) node[anchor=east]  {$\mathtt{D_n}$};
    \foreach \x in {-2,-1,0,...,3}
    \draw[thick,xshift=\x cm] (\x cm,-4) circle (1 mm);
    \foreach \y in {-2,-1,0,...,0,1}
    \draw[thick,xshift=\y cm] (\y cm,-4) ++(.3 cm, 0) -- +(14 mm,0);
		\draw[thick] (6.3,-3.9) -- +(1.4 cm, 0.5);
		\draw[thick] (6.3,-4.1) -- +(1.4 cm, -0.5);
    \draw[thick] (8 ,-3.3) circle (1mm);
		\draw[thick] (8 ,-4.7) circle (1mm);
		%\draw (7.8,0) node[anchor=east] {$\small{{\bf\boldsymbol{<}}}$};
		%\draw (-4,0) node[anchor=north]  {\tiny{1}};
		\draw (-4,-2.9) node[anchor=north]  {\tiny{2}};
   %\draw (-2,0) node[anchor=north]  {\tiny{2}};
		\draw (-2,-2.9) node[anchor=north]  {\tiny{2}};
    %\draw (0,0) node[anchor=north]  {\tiny{3}};
		\draw (0,-2.9) node[anchor=north]  {\tiny{4}};
    %\draw (2,0) node[anchor=north]  {\tiny{4}};
		\draw (2,-2.9) node[anchor=north]  {\tiny{4}};
    %\draw (4,0) node[anchor=north]  {\tiny{5}};
		\draw (4,-2.9) node[anchor=north]  {\tiny{6}};
    %\draw (5.9,1.5) node[anchor=north]  {\tiny{n-2}};
		\draw (5.8,-4.1) node[anchor=north]  {\tiny{{\tiny{$2\boldsymbol\lceil \frac{\textrm{n-2}}{2} \boldsymbol\rceil$}}}};
		%\draw (8.8,1) node[anchor=north]  {\tiny{n-1}};
    %\draw (7.55,-0.5) node[anchor=north]  {\tiny{n}};
		\draw (9.1,-2.4) node[anchor=north]  {\tiny{$\boldsymbol\lceil \frac{\textrm{n-1}}{2}\boldsymbol\rceil$}};
		\draw (9.1,-4.1) node[anchor=north]  {\tiny{$\boldsymbol\lceil \frac{\textrm{n-1}}{2}\boldsymbol\rceil$}};
    \draw[thick,dashed] (4.3,-4)--(5.7,-4);
    %\draw (4,3) node[anchor=north]  {\tiny{2}};
  %\end{tikzpicture}
%\end{center}
%%\bigskip
%\begin{center}
  %\begin{tikzpicture}[scale=.4]
    \draw (-21,-8) node[anchor=east]  {$\mathtt{E_6}$};
    \foreach \x in {-10,...,-6}
    \draw[thick,xshift=\x cm] (\x cm,-8) circle (1 mm);
    \foreach \y in {-10,...,-7}
    \draw[thick,xshift=\y cm] (\y cm,-8) ++(.3 cm, 0) -- +(14 mm,0);
    \draw[thick] (-16 cm, -6 cm) circle (1 mm);
    \draw[thick] ( -16 cm, -77mm) -- +(0, 1.4 cm);
    %\draw (0,0) node[anchor=north]  {\tiny{1}};
		\draw (-20,-6.9) node[anchor=north]  {\tiny{2}};
    %\draw (2,0) node[anchor=north]  {\tiny{3}};
		\draw (-18,-6.9) node[anchor=north]  {\tiny{4}};
    %\draw (4,-6.9) node[anchor=north]  {\tiny{4}};
		\draw (-16,-8.1) node[anchor=north]  {\tiny{6}};
    %\draw (6,0) node[anchor=north]  {\tiny{5}};
		\draw (-14,-6.9) node[anchor=north]  {\tiny{4}};
    %\draw (8,0) node[anchor=north]  {\tiny{6}};
		\draw (-12,-6.9) node[anchor=north]  {\tiny{2}};
    \draw (-15.4,-6) node[anchor=north]  {\tiny{2}};
		%\draw (4.4,2) node[anchor=north]  {\tiny{2}};
  %\end{tikzpicture}
%\end{center}
%\bigskip
%\begin{center}
  %\begin{tikzpicture}[scale=.4]
    \draw (-5,-8) node[anchor=east]  {$\mathtt{E_7}$};
    \foreach \x in {-2,...,3}
    \draw[thick,xshift=\x cm] (\x cm,-8) circle (1 mm);
    \foreach \y in {-2,...,2}
    \draw[thick,xshift=\y cm] (\y cm,-8) ++(.3 cm, 0) -- +(14 mm,0);
    \draw[thick] (0 cm,-6 cm) circle (1 mm);
    \draw[thick] (0 cm, -77mm) -- +(0, 1.4 cm);
    %\draw (0,0) node[anchor=north]  {\tiny{1}};
		\draw (-4,-6.9) node[anchor=north]  {\tiny{2}};
    %\draw (2,0) node[anchor=north]  {\tiny{3}};
    		\draw (-2,-6.9) node[anchor=north]  {\tiny{6}};
    %\draw (4,0) node[anchor=north]  {\tiny{4}};
    		\draw (0,-8.1) node[anchor=north]  {\tiny{8}};
    %\draw (6,0) node[anchor=north]  {\tiny{5}};
    		\draw (2,-6.9) node[anchor=north]  {\tiny{7}};
    %\draw (8,0) node[anchor=north]  {\tiny{6}};
    		\draw (4,-6.9) node[anchor=north]  {\tiny{4}};
   %\draw (10,0) node[anchor=north]  {\tiny{7}};
		\draw (6,-6.9) node[anchor=north]  {\tiny{3}};
    \draw (0.4,-6) node[anchor=north]  {\tiny{5}};
		%\draw (4.4,2) node[anchor=north]  {\tiny{2}};
  %\end{tikzpicture}
%\end{center}
%%\bigskip
%\begin{center}
  %\begin{tikzpicture}[scale=.4]
    \draw (-21,-12) node[anchor=east]  {$\mathtt{E_8}$};
    \foreach \x in {-10,...,-4}
    \draw[thick,xshift=\x cm] (\x cm,-12) circle (1 mm);
    \foreach \y in {-10,...,-5}
    \draw[thick,xshift=\y cm] (\y cm,-12) ++(.3 cm, 0) -- +(14 mm,0);
    \draw[thick] (-16 cm,-10 cm) circle (1 mm);
    \draw[thick] (-16 cm, -117mm) -- +(0, 1.4 cm);
    %\draw (0,0) node[anchor=north]  {\tiny{1}};
		\draw (-20,-10.9) node[anchor=north]  {\tiny{4}};
    %\draw (2,0) node[anchor=north]  {\tiny{3}};
    		\draw (-18,-10.9) node[anchor=north]  {\tiny{8}};
    %\draw (4,0) node[anchor=north]  {\tiny{4}};
    		\draw (-16,-12.1) node[anchor=north]  {\tiny{ $14$}};
    %\draw (6,0) node[anchor=north]  {\tiny{5}};
    		\draw (-14,-10.9) node[anchor=north]  {\tiny{$11$}};
    %\draw (8,0) node[anchor=north]  {\tiny{6}};
    		\draw (-12,-10.9) node[anchor=north]  {\tiny{$8$}};
   %\draw (10,0) node[anchor=north]  {\tiny{7}};
   		\draw (-10,-10.9) node[anchor=north]  {\tiny{6}};
   %\draw (12,0) node[anchor=north]  {\tiny{8}};
  		\draw (-8,-10.9) node[anchor=north]  {\tiny{2}};
    \draw (-15.5,-10) node[anchor=north]  {\tiny{$8$}};
		%\draw (4.5,2) node[anchor=north]  {\tiny{2}};
  %\end{tikzpicture}
%\end{center}
%\bigskip
%\begin{center}
  %\begin{tikzpicture}[scale=.4]
    \draw (-5,-12) node[anchor=east]  {$\mathtt{F_4}$};
		\draw[thick] (-4 ,-12) circle (1mm);
		\draw[thick] (-3.7cm,-12) -- +(1.4 cm, 0);
		\draw[thick] (-2 ,-12) circle (1mm);
    \draw[thick] (0,-12) circle (1mm);
    \draw[thick] (2 cm,-12) circle (1mm);
		\draw[thick] (0.3cm,-12) -- +(1.4 cm, 0);
    \draw[thick] (-17mm,-11.9) -- +(1.4 cm, 0);
    %\draw[thick] (0: 2.8 mm) -- +(0.7 cm, 0);
		\draw (-0.2,-12) node[anchor=east] {${{\bf\boldsymbol{>}}}$};
		%\draw[thick] (1mm) -- +(0.7 cm, 0);
    %\draw[thick] (0.9 cm, 0)[line width=0.5mm][<-] -- +(0.01 cm, 0);
    %\draw[thick] (0.9 cm, 0) -- +(0.88 cm, 0);
    \draw[thick] (-17mm,-12.1) -- +(1.4 cm, 0);
				%\draw (-2,0) node[anchor=north]  {\tiny{1}};
\draw (-4,-10.9) node[anchor=north]  {\tiny{2}};
		%\draw (0,0) node[anchor=north]  {\tiny{2}};
		\draw (-2,-10.9) node[anchor=north]  {\tiny{6}};
		%\draw (2,0) node[anchor=north]  {\tiny{3}};
		\draw (0,-10.9) node[anchor=north]  {\tiny{4}};
		%\draw (4,0) node[anchor=north]  {\tiny{4}};
		\draw (2,-10.9) node[anchor=north]  {\tiny{2}};
  %\end{tikzpicture}
%\end{center}
%\bigskip
%\begin{center}
  %\begin{tikzpicture}[scale=.4]
    \draw (5,-12) node[anchor=east]  {$\mathtt{G_2}$};
    \draw[thick] (6 ,-12) circle (1mm);
    \draw[thick] (8 cm,-12) circle (1mm);
    \draw[thick] (6.3,-11.9) -- +(1.4 cm, 0);
    %\draw[thick] (0: 2.8 mm) -- +(0.7 cm, 0);
		\draw (7.8,-12) node[anchor=east] {${{\bf\boldsymbol{<}}}$};
		%\draw[thick] (1mm) -- +(0.7 cm, 0);
    %\draw[thick] (0.9 cm, 0)[line width=0.5mm][<-] -- +(0.01 cm, 0);
    %\draw[thick] (0.9 cm, 0) -- +(0.88 cm, 0);
    \draw[thick] (6.3,-12) -- +(1.4 cm, 0);
		\draw[thick] (6.3,-12.1) -- +(1.4 cm, 0);
		%\draw (0,0) node[anchor=north]  {\tiny{1}};
		\draw (6,-10.9) node[anchor=north]  {\tiny{2}};
		%\draw (2,0) node[anchor=north]  {\tiny{2}};
		\draw (8,-10.9) node[anchor=north]  {\tiny{2}};
  \end{tikzpicture}
\end{center}
\end{theom}

The paper is organized as follows: In Section~\ref{first} we introduce definitions and basic notations, in Section~\ref{three} we prove Theorem~\ref{B}.\\[4mm]
\textbf{Acknowledgements}

T.B. was funded by the DFG-priority program 1388 “Representation Theory”, grant ''LI 990/10-1", G.F. was partially funded the grant "FO 867/1-1",  L.B. and  C.D. were partially funded within the framework of this program. The main work of this article has been conducted during a workshop organized at the University of Glasgow, and all authors would like to thank the Glasgow Mathematical Journal Trust Fund, the Edinburgh Mathematical Journal and especially the University of Glasgow for this opportunity.\\[4mm]
Mathematisches Institut, Universit\"at zu K\"oln, Germany (T.B., L.B., C.D.)\\
%\email{tbackha@math.uni-koeln.de}
%\email{lbossing@math.uni-koeln.de}
%\email{cdesczyk@math.uni-koeln.de}
%\address{\newline
Mathematisches Institut, Universit\"at Bonn, Germany, School of Mathematics and Statistics, University of Glasgow, UK (G.F.)
%School of Mathematics and Statistics, University of Glasgow, UK}

%%%%%%%%%%%%%%%%%%%%%%%%%%%%%%%%%%%%%%%%%%%%%%%%%%%%%%%%%%%%%%%%%%%%%%%%%%%%%%%%%%%%%%%%%%%%%%%%%%%%%%%%%%%%%%%%%%%%%%%%%%%%%%%%%%%%
%%%%%%%%%%%%%%%%%%%%%%%%%%%%%%%%%%%%%%%%%%%%%%%%%%%%%%%%%%%%%%%%%%%%%%%%%%%%%%%%%%%%%%%%%%%%%%%%%%%%%%%%%%%%%%%%%%%%%%%%%%%%%%%%%

\section{Preliminairies}\label{first}
Let $\mathfrak g$ be a simple Lie algebra of rank $n$. We fix a Cartan subalgebra $\Lh$ and a triangular decomposition $\mathfrak g=\mathfrak n^+ \oplus \mathfrak h \oplus \mathfrak n^-$.
The set of roots (resp. positive roots) of $\Lg$ is denoted  $R$ (resp. $R^+$), $\theta$ denotes the highest root. Let $\alpha_i, \om_i \ i=1,...,n$ be the simple roots and the fundamental weights. Let $W$ be the Weyl group associated to the simple roots and $w_0\in W$ the longest element. For $\alpha \in R^+$ we fix a $\mathfrak{sl}_2$ triple $\{ e_\alpha, f_{\alpha}, h_{\alpha}= [e_\alpha, f_{\alpha}]\}$. The integral weights and the dominant integral weights are denoted $P$ and $P^+$. \\
Let $\lbrace x_1,x_2,...\rbrace$ be an ordered basis of $\Lg$, then $U(\Lg)$ denotes the universal enveloping algebra of $\Lg$ with PBW basis $\lbrace x_{i_1} \cdot \cdot \cdot x_{i_m}\mid m\in \mathbb Z_{\ge 0}, i_1\le i_2\le ... \le i_m\rbrace$.

 \subsection{Modules}
 For $\lambda \in P^+$ we consider the irreducible $\mathfrak g$-Module $V(\lam)$ with highest weight $\lambda$. Then $V(\lambda)$ admits a decomposition into $\Lh$-weight spaces, \\ ${V(\lam)=\bigoplus_{\tau\in P}V(\lam)_\tau}$ with $V(\lam)_{\lam}$ and $V(\lam)_{w_0(\lam)}$, the highest and lowest weight spaces, being one dimensional. Let $v_\lam$ denote the highest weight vector, $v_{w_0(\lam)}$ denote the lowest weight vector satisfying
\[
e_\alpha v_\lam=0, \; \, \forall \, \alpha \in R^+ \; ; \; f_\alpha v_{w_0(\lam)}=0, \, \, \forall \, \alpha \in R^+.
\]
We have  $U(\Ln ^-).v_\lam\cong V(\lam)\cong U(\lie n^+).v_{w_0(\lam)}$.\\
The comultiplication $(x \mapsto x \otimes 1 + 1 \otimes x)$ provides a $\lie g$-module structure on $V(\lambda) \otimes V(\mu)$. This module decomposes into irreducible components, where the Cartan component generated by the highest weight vector $v_\lambda \otimes v_\mu$ is isomorphic to $V(\lambda + \mu)$.

 \subsection{PBW-filtration}
The Hilbert-Poincar\'e series of the PBW-graded module $V(\lam)^a:=\bigoplus_{s\ge 0} V(\lam)_s /V(\lam)_{s-1}$ is the polynomial
\begin{eqnarray*}
p_\lam(q)&=&\sum\nolimits_{s\ge0}  \dim(V(\lam)_s/V(\lam)_{s-1})q^s\\
&=&1+\dim(V(\lam)_1/V(\lam)_0)q+\dim(V(\lam)_2/V(\lam)_1)q^2+...\nonumber
\end{eqnarray*}
and we define the PBW-degree of $V(\lam)$ to be $\deg(p_\lam(q))$.\\

It is easy to see that  $\Ln^+.( U(\Ln^-)_s.v_\lam)\subseteq U(\Ln^-)_s.v_\lam \ \forall \, s\ge 0$ (see also \cite{FFoL11a})  and hence
%\begin{eqnarray}\label{n-stable}
$U(\Ln^+).V(\lam)_s\subseteq V(\lam)_s.$
%\end{eqnarray}
Let $s_\lam$ be minimal such that $v_{w_0(\lam)}\in V(\lam)_{s_\lam}$. Then ${V(\lam)=U(\Ln^+).v_{w_0(\lam)}\subseteq V(\lam)_{s_\lam}}$ and
\begin{coro}\label{s-lambda} $s_\lambda = \operatorname{deg} (p_\lambda(q))$ and
\begin{eqnarray*}
V(\lam)=V(\lam)_{s_\lam}.
\end{eqnarray*}
\end{coro}

\subsection{Graded weight spaces}
The PBW filtration is compatible with the decomposition into $\lie h$-weight spaces:
\[
\dim V(\lambda)_\tau = \sum_{s\ge0} \dim \left(V(\lambda)_s/V(\lambda)_{s-1}\right) \cap V(\lambda)_\tau.
\]
So we can define for every weight $\tau$ the Hilbert-Poincar\'e polynomial:
\[
p_{\lambda, \tau}(q) = \sum_{s \geq 0} \dim \left(V(\lambda)_s/V(\lambda)_{s-1}\right)_{\tau}q^s \text{ and then }
p_{\lambda}(q) = \sum_{\tau \in P } p_{\lambda,\tau}(q).
\]
A natural question is, if we can extend our results to these polynomials? If the weight space $V(\lambda)_\tau$ is one-dimensional, then $p_{\lambda,\tau}(q)$ is a power of $q$. For $\tau = \lambda$ this is constant $1$, for $\tau = w_0(\lambda)$, the lowest weight, this is $q^{\deg  p_{\lambda}(q) }$ as we have seen in  Corollary~\ref{s-lambda}. A first approach to study these polynomials can be found in \cite{CF13}.\\

\subsection{Graded Kostant partition function}
For the readers convienience we recall 
%the Kostant partition function $P(\nu)$ (see \cite{Ko59}), which counts the number of decompositions of a fixed weight into a sum of positive roots, and how it is related to our study. For this, we consider the power series
%\[
%\prod_{\alpha > 0} \frac{1}{(1 - e^{\alpha} )}
%\]
%and its expansion
%\[
%\sum_{\nu\in P} P_{\nu} e^{\nu}.
%\]
%\\
%We have immediately that  $P_{\nu}$ is the dimension of the weight space of weight $- \nu$ in the $\lie h$-module  $U(\lie n^-)$. Since we are interested in graded modules, we consider 
here the \textit{graded Kostant partition function} (see \cite{Kos59}), which counts the number of decompositions of a fixed weight into a sum of positive roots, and how it is related to our study. We consider the power series and its expansion:
\[
\prod_{\alpha > 0} \frac{1}{(1 - qe^{\alpha} )}, \,\,\,\,\,\,
%\]
%\math{and its expansion}
%\[
\sum_{\nu\in P} P_{\nu}(q) e^{\nu}.
\]
We have immediately $\chara S(\lie n^-) = \sum_{\nu \in P} P_{\nu}(q) e^{-\nu}$.
\begin{remi}
For a polynomial $p(q) = \sum_{i = 0}^{n} a_i q^i$, we denote $\mindeg p(q)$ the minimal $j$ such that $a_j \neq 0$. Then we have obviously
\begin{equation}
\mindeg p_{\lambda, \nu}(q) \geq \mindeg P_{\lambda - \nu}(q).
\label{eq:min-deg}
\end{equation}
\end{remi}
We will use this inequality for the very special case $\nu = w_0(\lambda)$ in the proof of Theorem~\ref{B}.\\
We see from Theorem \ref{B} that this inequality is a proper inequality for certain cases in exceptional type as well as $B_n D_n$ (this has been noticed also in \cite{CF13}).

%%%%%%%%%%%%%%%%%%%%%%%%%%%%%%%%%%%%%%%%%%%%%%%%%%%%%%%%%%%%%%%%%%%%%%%%%%%%%%%%%%%%%%%%%%%%%%%%%%%%%%%%%%%%%%%%%%%%%%%%%%%%%%%%%%%
%         Dynkin
%%%%%%%%%%%%%%%%%%%%%%%%%%%%%%%%%%%%%%%%%%%%%%%%%%%%%%%%%%%%%%%%%%%%%%%%%%%%%%%%%%%%%%%%%%%%%%%%%%%%%%%%%%%%%%%%%%%%%%%%%%%%%%%%%%%
\section{Proof of Theorem \ref{B}}\label{three}
In this section we will provide a proof of Theorem~\ref{B}. For a fixed $1 \leq i \leq \operatorname{ rank } \lie g$, we will give a monomial $u \in U(\lie n^-)$ of the predicted degree mapping the highest weight vector $v_{\omega_i}$ to the lowest weight vector $v_{w_0(\omega_i)}$. We then show that there is no monomial of smaller degree satisfying this.\\ To write down these monomials explicitly, let us denote $\theta_{X_n}$ the highest root of a Lie algebra of type $X_n$. We set further (using the indexing from \cite{Hum72}):
\begin{itemize}
\item In the $A_n$-case, $Y_{n-2}$ the type of the Lie algebra generated by the simple roots $\{ \alpha_{2}, \ldots ,\alpha_{n-1} \}$.
\item In the $B_n, D_n$-case, $Y_{n-k}$ the type of the Lie algebra generated by the simple roots $\{ \alpha_{k+1}, \ldots ,\alpha_{n} \}$.
\item In the exceptional and symplectic cases, $\theta_{X_n} = c_k \omega_k$ for some $k$, $Y_{n-1}$ the type of the Lie algebra generated by the simple roots $\{\alpha_1, \ldots, \alpha_n\}\setminus\{ \alpha_k\}$.
\end{itemize}
Let $u \in U(\lie n^-)$ be one of the monomials in Figure \ref{monoms}. It can be seen easily from Figure \ref{monoms} that $u = f_{\theta_{X_n}}^{a_i^{\vee}}u_1$, where $a_i^{\vee} = w_i(h_{\theta_{X_n}})$ and $u_1$ is the monomial in Figure \ref{monoms} corresponding to the restriction of $\omega_i$ to the Lie subalgebra of type  $Y_{n-\ell}$. If we denote $\lie n_1^-$ the lower part in the triangular decomposition of the Lie subalgebra of type $Y_{n-\ell}$, then $u_1 \in U(\lie n^-_1)$.\\ 
Let $u=f_{\theta_1}^{b_1}f_{\theta_2}^{b_2}\dots f_{\theta_r}^{b_r}$. Note that all $f_{\theta_j}$ commute and it is easy to see that $\theta_j(h_{\theta_{j+p}})=0,\ \forall p \geq 0$ (since $\theta_j$ is a sum of fundamental weights, which are all orthogonal to the simple roots of the Lie algebra with highest root $\theta_{j+p}$) and $b_j = \omega_i (h_{\theta_j})$. \\
The Weyl group $W$ acts on $V(\omega_i)$ and if $v$ is an extremal weight vector of weight $\mu$, then $w.v$ is a nonzero extremal weight vector of weight $w(\mu)$. Further if $w=s_\alpha$ (reflection at a root $\alpha$) and $\mu(h_\alpha) \geq 0$, then $w.v=c^*f_{\alpha}^{\mu(h_\alpha)}.v$ for some $c^* \in \mathbb C^*$.\\
Now consider $w=s_{\theta_r} \dots s_{\theta_1}$, where $s_{\theta_j}$ is the reflection at the root $\theta_j$. Then we have $w.v_{\omega_i}=v_{w_0(\omega_i)}=u.v_{\omega_i} \neq 0$ in $V(\omega_i)$.
\begin{figure}
$$
\begin{array}{lll}
X_n & \om_i = \theta_{X_n} & f_{\theta_{X_n}}^2\\
A_n &  \om_i & f_{\theta_{A_n}}f_{\theta_{A_{n-2}}} \cdots f_{\theta_{A_{n + 2  - 2 \operatorname{min} \{i, n-i\}}}}\\
C_n & \om_i & f_{\theta_{C_n}}f_{\theta_{C_{n-1}}} \cdots f_{\theta_{C_{n + 1 -i}}}\\
B_n & \om_{2i} &  f^2_{\theta_{B_n}}f^2_{\theta_{B_{n-2}}} \cdots f^2_{\theta_{B_{n + 2 -2i}}}\\
B_n  &\om_{2i+1}& f^2_{\theta_{B_n}}f^2_{\theta_{B_{n-2}}} \cdots f^2_{\theta_{B_{n -2i}}} f_{\alpha_{2i+1}}\\
B_n & n \text{ even} , \om_n & f_{\theta_{B_n}} f_{\theta_{B_{n-2}}} \cdots f_{\theta_{B_2}}\\
B_n & n \text{ odd} , \om_n & f_{\theta_{B_n}}f_{\theta_{B_{n-2}}}  \cdots f_{\theta_{B_2}} f_{\alpha_n}\\
D_n & \om_{2i} &f^2_{\theta_{D_n}}f^2_{\theta_{D_{n-2}}} \cdots f^2_{\theta_{D_{n + 2 -2i}}}\\
D_n & \om_{2i+1} &  f^2_{\theta_{D_n}}f^2_{\theta_{D_{n-2}}} \cdots f_{\theta_{D_{n - 2i}}} f_{\alpha_{2i+1}}\\
D_n & n \text { even}, \om_i, i = n-1, n &  f_{\theta_{D_n}} f_{\theta_{D_{n-2}}} \cdots f_{\theta_{D_4}} f_{\alpha_{i}}\\
D_n & n \text { odd}, \om_i, i = n-1, n &  f_{\theta_{D_n}} f_{\theta_{D_{n-2}}} \cdots f_{\theta_{D_{5}}}f_{\theta_{A_4}}\\
E_6 & \om_1, \om_6 & f_{\theta_{E_6}} f_{\theta_{A_5}}\\
E_6 & \om_3, \om_5 & f_{\theta_{E_6}}^2 f_{\theta_{A_5}} f_{\theta_{A_3}}\\
E_6 & \om_4 & f_{\theta_{E_6}}^3f_{\theta_{A_5}}f_{\theta_{A_3}}f_{\alpha_4}\\
%E_7 & \om_1 & f_{\theta_{E_7}}^2\\
E_7 & \om_2 &  f_{\theta_{E_7}}^2 f_{\theta_{D_6}} f_{\theta_{D_4}}f_{\alpha_2}\\
E_7 & \om_3 &   f_{\theta_{E_7}}^3 f_{\theta_{D_6}} f_{\theta_{D_4}}f_{\alpha_3}\\
E_7 & \om_4 & f_{\theta_{E_7}}^4 f_{\theta_{D_6}}^2 f_{\theta_{D_4}}^2\\
E_7 & \om_5 &  f_{\theta_{E_7}}^3 f_{\theta_{D_6}}^2 f_{\theta_{D_4}} f_{\alpha_5}\\
E_7 & \om_6 &  f_{\theta_{E_7}}^2 f_{\theta_{D_6}}^2\\
E_7 & \om_7 & f_{\theta_{E_7}} f_{\theta_{D_6}} f_{\alpha_7}\\
E_8 & \om_1 &   f_{\theta_{E_8}}^2 f_{\theta_{E_7}}^2\\
E_8 & \om_2 &  f_{\theta_{E_8}}^3 f_{\theta_{E_7}}^2 f_{\theta_{D_6}} f_{\theta_{D_4}}f_{\alpha_2}\\
E_8 & \om_3 &   f_{\theta_{E_8}}^4 f_{\theta_{E_7}}^3 f_{\theta_{D_6}} f_{\theta_{D_4}}f_{\alpha_3}\\
E_8 & \om_4 & f_{\theta_{E_8}}^6  f_{\theta_{E_7}}^4 f_{\theta_{D_6}}^2 f_{\theta_{D_4}}^2\\
E_8 & \om_5 &  f_{\theta_{E_8}}^5 f_{\theta_{E_7}}^3 f_{\theta_{D_6}}^2 f^{}_{\theta_{D_4}} f_{\alpha_5} \\
E_8 & \om_6 &  f_{\theta_{E_8}}^4 f_{\theta_{E_7}}^2 f_{\theta_{D_6}}^2\\
E_8 & \om_7 & f_{\theta_{E_8}}^3f_{\theta_{E_7}} f_{\theta_{D_6}} f_{\alpha_7} \\
%E_8 & \om_8 & f_{\theta_{E_8}}^2 \\
%F_4 & \om_1 & f_{\theta_{F_4}}^2\\
F_4 & \om_2 & f_{\theta_{F_4}}^3 f_{\theta_{C_3}} f_{\theta_{A_2}} f_{\alpha_2}\\
F_4 & \om_3 & f_{\theta_{F_4}}^2 f_{\theta_{C_3}}f_{\theta_{C_2}}\\
F_4 & \om_4 & f_{\theta_{F_4}}f_{\theta_{C_3}}\\
G_2 & \om_1 & f_{\theta_{G_2}}f_{\alpha_1}
%G_2 & \om_2 & f_{\theta_{G_2}}^2
\end{array}
$$
\caption{}
\label{monoms}
\end{figure}
So we obtain an upper estimate for the degree.\\
In general the degree of $u$ is bigger than the minimal degree coming from Kostant's graded partition function \eqref{eq:min-deg}. For $A_n, C_n$ the degrees coincide and hence we are done in these cases.\\
We will prove Theorem~\ref{B} for the remaining cases $X_n$ by induction on the rank of the Lie algebra. So we want to prove that if $p \in U(\mathfrak{n}^-)$ with $p.v_{\omega_i} = v_{w_0(\omega_i)}$ then $\deg(p) \geq \deg(u)$, where $u$ is from Figure 1.\\
Consider the induction start, e.g. $\omega_i = \theta_{X_n}$, then the minimal degree is obviously $2$. The maximal non-vanishing power of $f_{\theta_{X_n}}$ is certainly $a_i^{\vee}$
%, where $h_\theta = \sum a_j^{\vee} h_j$,
and $f_{\theta_{X_n}}^{a_i^{\vee}}.v_{\omega_i}$ is the highest weight vector of a simple module of fundamental weight for the Lie algebra $Y_{n-l}$ defined as above. By induction we know that if $q \in U(\mathfrak{n}_1^-)$ with $q.(f_{\theta_{X_n}}^{a_i^{\vee}}.v_{\omega_i}) = v_{w_0(\omega_i)}$ then $\deg(q) \geq \deg(u_1)$.\\
First we suppose $f_{\theta_{X_n}}^{a_i^{\vee}}.v_{\omega_i}$ is a factor of $p$, so $p = f_{\theta_{X_n}}^{a_i^{\vee}}p'$ and then by weight considerations $p' \in U(\mathfrak{n}^-_1)$. Then $p'.(f_{\theta_{X_n}}^{a_i^{\vee}}.v_{\omega_i}) = v_{w_0(\omega_i)}$ (the lowest weight vector in $V(\omega_i)$ as well as in the simple submodule). Therefore $\deg(p') \geq \deg(u_1)$ which implies $\deg(p) \geq \deg(u)$.\\
Suppose now the maximal power of $f_{\theta_{X_n}}$ in $p$ is $f_{\theta_{X_n}}^{a_i^{\vee}-k}$, $k \geq 0$ and $\deg(p) < \deg(u)$. Let $X_n$ be of type $B_n, D_n$ or exceptional, then $\theta_{X_n} = \om_j$ and we denote 
\[
R^+_s = \{ \alpha \in R^+ \, | \, w_j(h_{\al}) = s\},
\]
Then $R_2^+ = \{ \theta_{X_{n}} \}$ and if $\beta\in R^+_1$ then $\theta_{X_n} - \beta \in R^+_1$. By weight reasons $p = f_{\theta_{X_n}}^{a_i^{\vee} - k} f_{\beta_1} \cdots f_{\beta_{2k}} p_1$ for some $\beta_1, \ldots, \beta_{2k} \in R^+_1$ and some polynomial $p_1$ in root vectors of roots in $R_0^+$. We have to show that  $p.v_{\om_i} = 0 \in V(\om_i)^a$ and we will use induction on $k$ for that:
The induction start is $k = 0$. The induction step is for $k \geq 1$:
$$
\begin{array}{rcl}
0 = p_1 f_{\theta_{X_n}}^{a_i^{\vee} + k}.v_{\om_i} & = & (e_{\theta_{X_n} - \beta_1}) \cdots (e_{\theta_{X_n} - \beta_{2k}}) p_1 f_{\theta_{X_n}}^{a_i^{\vee} + k}.v_{\om_i} \\
& = &  cf_{\theta_{X_n}}^{a_i^{\vee} - k} f_{\beta_1} \cdots f_{\beta_{2k}} p_1 .v_{\om_i} + \sum_{\ell > 0}^k   f_{\theta_{X_n}}^{a_i^{\vee} - k + \ell} q_{\ell} .v_{\om_i}
\end{array}
$$
for some $c \in \mathbb{C}^*, q_{\ell} \in U(\lie n^-)$. For $0 \leq \ell < k$ all the summands are equals to zero by induction (on $k$). For $\ell = k$, we recall our assumption  $\deg (p) < \deg (u)$ and so $\deg (q_k) < \deg (u_1)$ which implies $f^{a_i^\vee}q_k.v_{\omega_i} = 0$. So we can conclude
$f_{\theta_{X_n}}^{a_i^{\vee} - k} f_{\beta_1} \cdots f_{\beta_{2k}} p_1 .v_{\om_i}=0$.
\bibliographystyle{alpha}
\bibliography{../pbw-demazure-bib}

\end{document}